\documentclass[]{amsart}

\usepackage{amsmath, amssymb, amsthm, amsfonts}

\newtheorem{theorem}{Theorem}[section]
\newtheorem{lemma}[theorem]{Lemma}

\theoremstyle{definition}
\newtheorem{definition}[theorem]{Definition}
\newtheorem{remark}[theorem]{Remark}
\newtheorem*{ack}{Acknowledgements}

\usepackage{graphicx}

\usepackage{overpic}

\begin{document}

\title{Genus zero PALF structures on the Akbulut-Yasui plugs}
\author{Takuya Ukida}
\address{Department of Mathematics, Tokyo Institute of Technology,
2-12-1 Oh-okayama, Meguro-ku, Tokyo 152-8551, Japan}
\date{}

\begin{abstract}
We construct a genus zero PALF structure 
on each of plugs introduced by Akbulut and Yasui 
and describe the monodromy as a positive factorization in the mapping class group 
of a fiber. 
We also examine the monodromies of PALFs on 
a certain pair of compact Stein surfaces such that one is obtained 
by applying a plug twist to the other. 
\end{abstract}

\maketitle

\section{Introduction}

The problem of classifying all differential structures 
defined on a given 4-manifold 
is an important problem 
in understanding the overall picture of a 4-manifold.

Akbulut and Yasui \cite{AY1} introduced corks and plugs.
Corks and plugs are compact Stein surfaces. 
Matveyev, Curtis-Freedman-Hsiang-Stong, and Akbulut-Matveyev's theorem show 
that the study of exotic manifold pairs constructed using cork 
is important for the classification problem 
of the differential structure of 4-manifolds. 

\begin{theorem}(Matveyev \cite{Matveyev}, Curtis-Freedman-Hsiang-Stong \cite{CFHS}, Akbulut-Matveyev \cite{AM})
    For every homeomorphic but non-diffeomorphic pair of simply connected
    closed 4-manifolds, 
    one is obtained from the other by removing a contractible 
    4-manifold and gluing it via an involution on the boundary. 
    Such a contractible 4-manifold has since been called a Cork. 
    Furthermore, corks and their complements can always
    be made compact Stein 4-manifolds.
\end{theorem}

It is shown by Akbulut and Yasui \cite{AY3} using cork
that an infinite number of exotic Stein surface pairs 
embedded in X exist for any four-dimensional two-handle body X with $b_2(X) \geq 1$.

The plug generalizes the Gluck twist.
The plug is also used to make exotic manifolds, as well as cork.

On the other hand, 
Loi and Piergallini \cite{LP} proved that every compact Stein surface admits
a positive allowable Lefschetz fibration over $D^2$ (a PALF for short).
Therefore we can investigate compact 
Stein surfaces in terms of positive factorizations in mapping class groups 
(see also Akbulut and Ozbagci \cite{AO}, Akbulut and Arikan \cite{AA}).

Since corks and plugs are Stein surfaces, 
the study of the relationship
between Stein surfaces and mapping class groups using PALFs
plays an important role
in classifying differential structures on 4-manifolds.

If a PALF is created from a Stein surface by the existing method ( \cite{LP}, \cite{AO}, \cite{AA} ),
its genus will be large,
and it will be complicated and difficult
to handle as a mapping class group element.

Gompf \cite{Go} indicates that the Stein surface is compatible with Kirby calculus.
In this paper, we use Kirby calculus to make the best PALF with genus zero from the Akbulut-Yasui plug.

One planar (i.e. genus zero) PALF on the Akbulut cork was made in the previous paper of the author \cite{Uki1}, 
but in this paper, we made an infinite number of planar PALFs on the Akbulut-Yasui plugs. 
Being planar is also playing an important role in \cite{KOU}.

In this paper,  we first construct a genus zero PALF structure 
on each of plugs introduced by Akbulut and Yasui \cite{AY1} 
and describe the monodromy as a positive factorization in the mapping class group 
of a fiber.

\begin{figure}[htbp]
 \begin{center}
  \includegraphics[width=75mm]{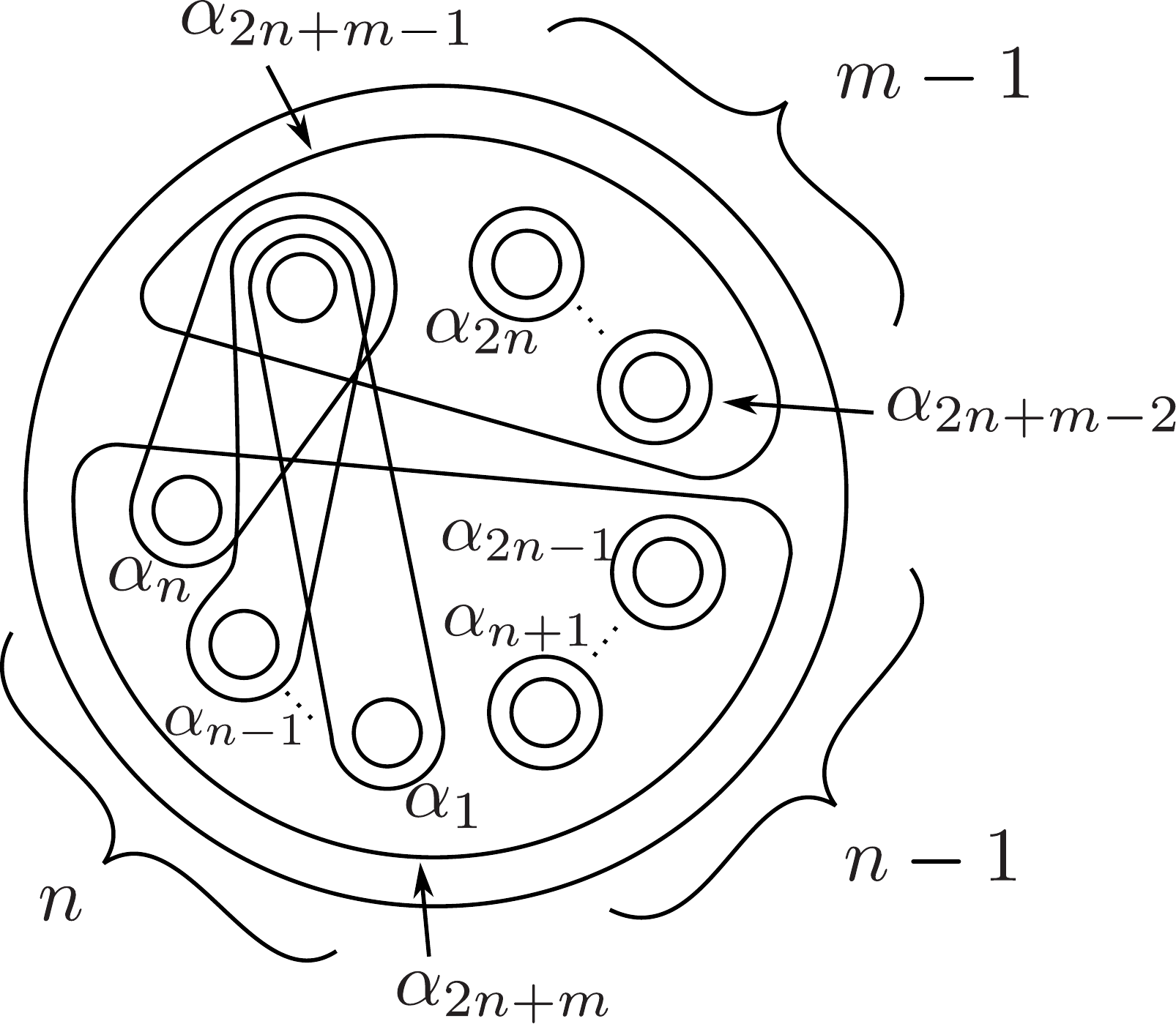}
 \end{center}
    \caption{Vanishing cycles of a genus zero PALF on $W_{m, n}$.}
 \label{AkYsPlug}
\end{figure}

\begin{theorem}\label{MainThm:plug1}
    For any $m \ge 1, n \ge 2$,
    the Akbulut-Yasui plug $(W_{m,n}, f_{m,n})$
    admits a genus zero PALF structure.
    The monodromy of the  PALF is described by the factorization
    $t_{\alpha_{2n+m}} \cdots t_{\alpha_1}$, 
    where $t_\alpha$ is a right-handed Dehn twist along a
    simple closed curve $\alpha$ on a fiber and
    $\alpha_{2n+m}, \ldots, \alpha_1$ are simple closed curves
    shown in Figure \ref{AkYsPlug}.
\end{theorem}

\begin{figure}[htbp]
    \begin{center}
        \includegraphics[width=30mm]{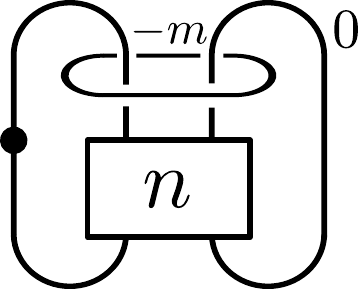}
    \end{center}
    \caption{Kirby diagram for $W_{m, n}$}
    \label{AY_plug_KirbyDiam}
\end{figure}

\begin{figure}[htbp]
 \begin{center}
  \includegraphics[width=85mm]{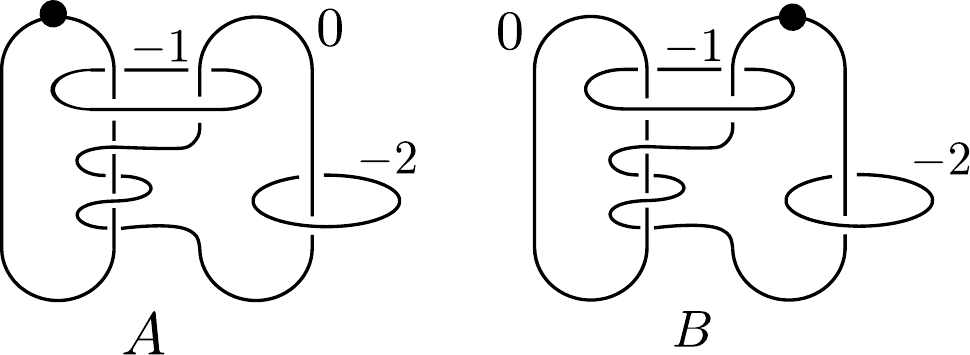}
 \end{center}
 \caption{Kirby diagrams for $A$ and $B$.}
 \label{MfdAB}
\end{figure}

Note that the genus of a PALF on the manifold $W_{m, n}$ 
in a known way (cf. \cite{AO} and  \cite{AA}) is much more than zero.
We obtained similar results for the Akbulut cork $W_1$ \cite{Uki1}.
In the present paper, we construct a genus zero PALF on an infinite number of 
the Akbulut-Yasui plugs.

\begin{figure}[htbp]
 \begin{center}
  \includegraphics[width=45mm]{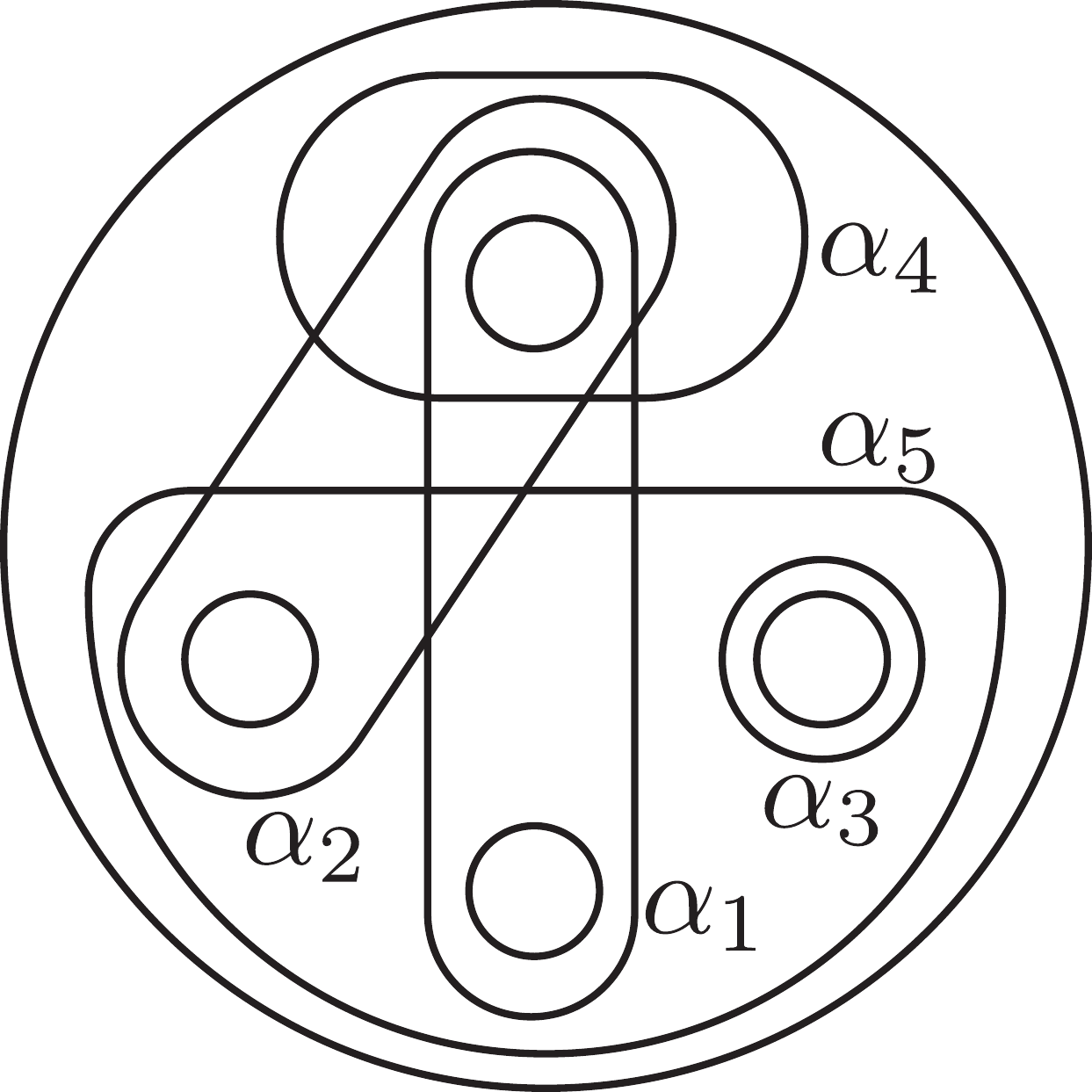}
 \end{center}
 \caption{Vanishing cycles of a genus zero PALF on $W_{1, 2}$.}
 \label{W12}
\end{figure}

\begin{figure}[]
    \begin{center}
        \includegraphics[width=45mm]{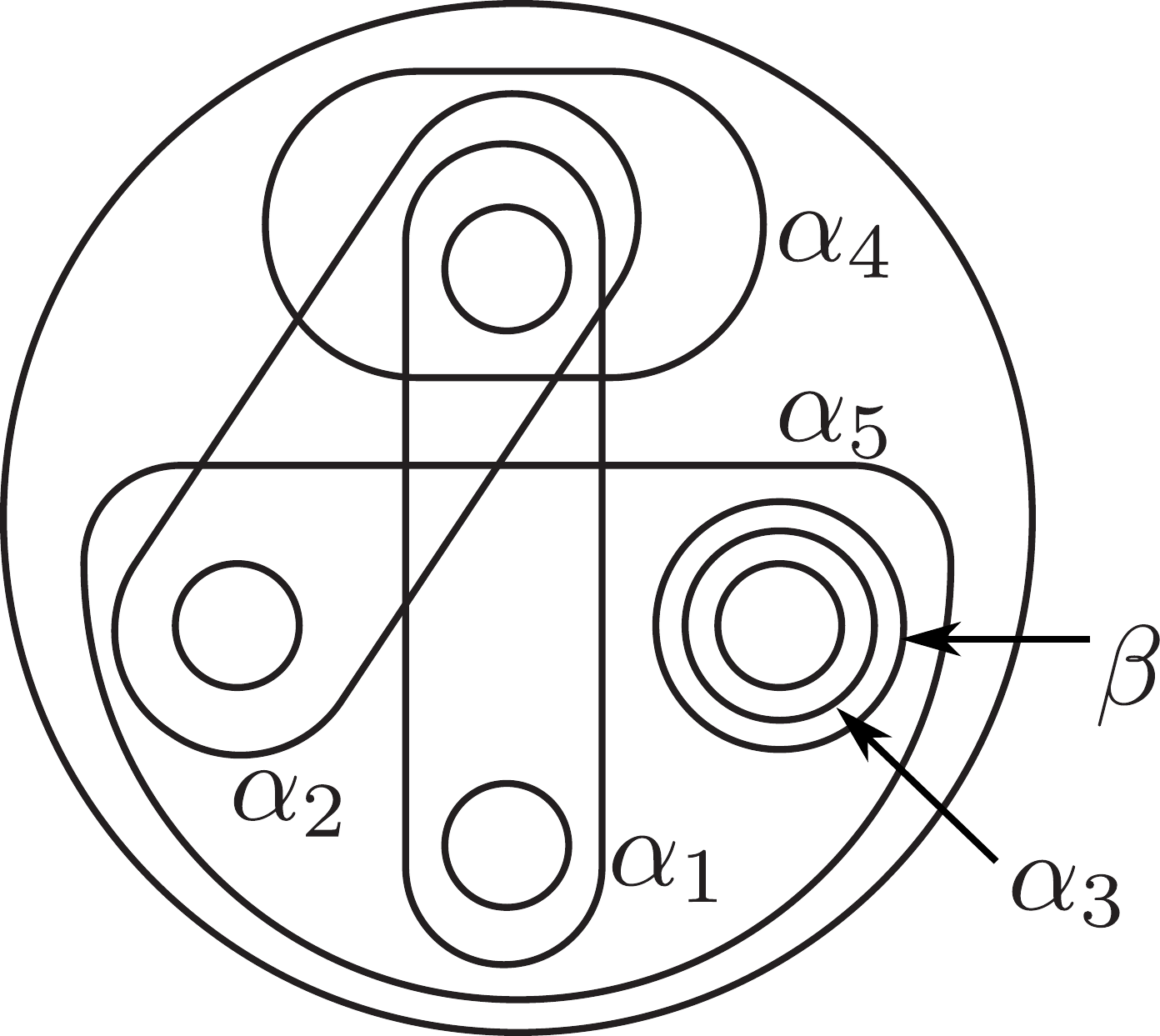}
    \end{center}
    \caption{Vanishing cycles of a genus zero PALF.}
    \label{Thm1-2ProofPALF_A}
\end{figure}

\begin{figure}[]
    \begin{center}
        \includegraphics[width=45mm]{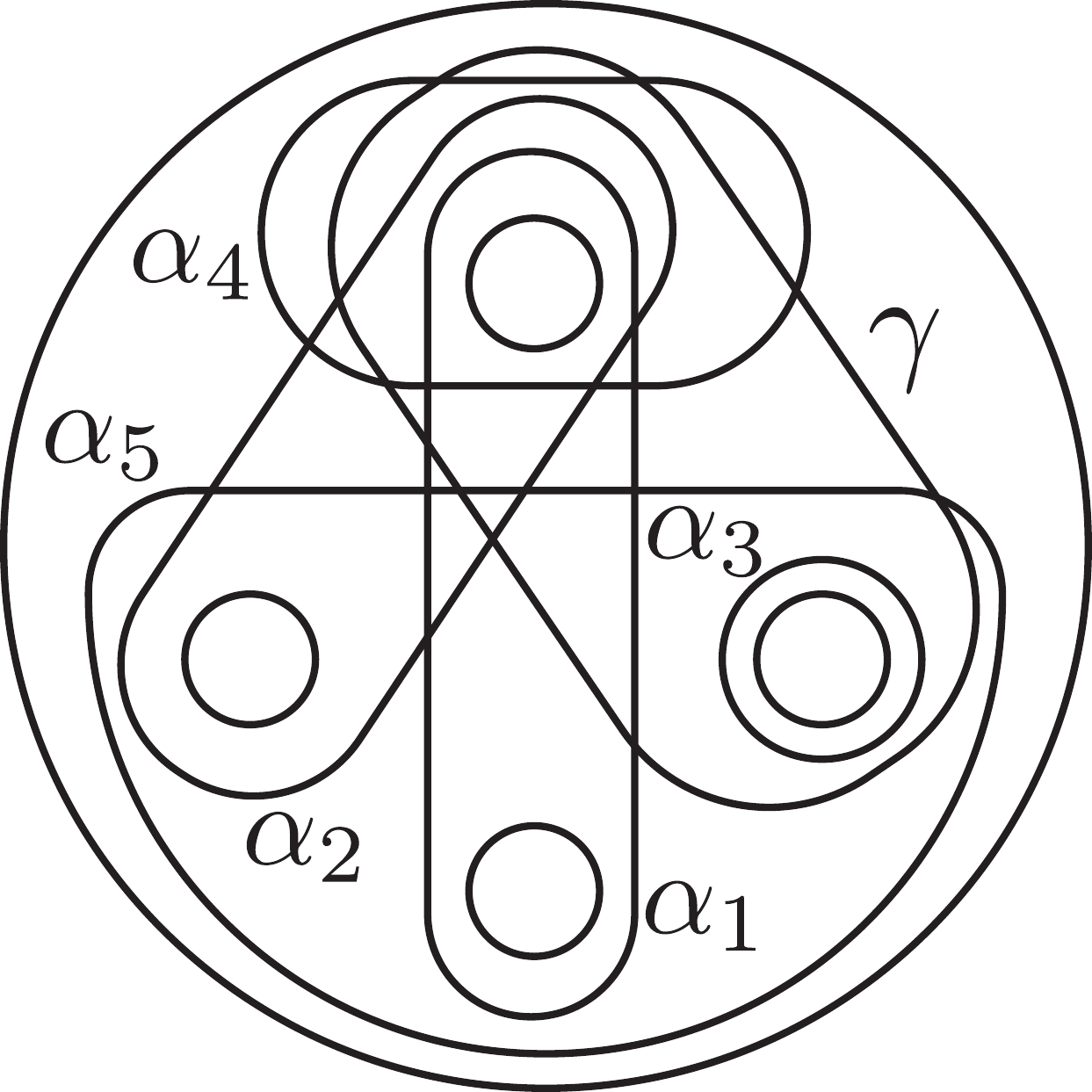}
    \end{center}
    \caption{Vanishing cycles of a genus one PALF.}
    \label{Thm1-2ProofPALF_B}
\end{figure}

Furthermore, we show that
example of two 4-manifolds $A$ and $B$ which is obtained from $A$ by
plug twist of $A$.
The manifolds $A$ and $B$ admit genus zero PALF structure, and have the following properties:

\begin{theorem}\label{MainThm:plug2}
    The manifolds $A$ and $B$ which are showed by the Kirby diagrams in the
    Figure~\ref{MfdAB} have the following properties.
    \begin{enumerate}
        \renewcommand{\labelenumi}{(\arabic{enumi})}
        \setlength{\itemsep}{0cm} 
        \setlength{\parskip}{1mm} 
        \item $B$ is obtained by  plug twist of $A$ along
              the Akbulut-Yasui plug $(W_{1, 2}, f_{1, 2})$.
        \item $A, B$ admit genus zero PALF structure.
        \item The second Betti numbers of $A$ and $B$ are $2$, and the second homology groups of
        $A$ and $B$ are isomorphic.
        \item The boundaries of $A$ and $B$ are diffeomorphic.
        \item $A, B$ do not have isomorphic intersection numbers. 
              Especially $A$ and $B$ are not homeomorphic.
        \item The monodromy representation of
	      genus zero PALF structures which $W_{1, 2}$ admits
              is $t_{\alpha_4} t_{\alpha_3} t_{\alpha_2} t_{\alpha_1}$, 
              where $\alpha_i$ is a simple closed curve in diagram \ref{W12},
              and $t_{\alpha_i}$ is right-handed Dehn twist along $\alpha_i$.
              Monodromy representations of
              genus zero PALF structure of $A, B$ are 
              $t_{\alpha_4} t_{\alpha_3} t_{\alpha_2} t_{\alpha_1} t_{\beta}$ and
              $t_{\alpha_4} t_{\alpha_3} t_{\alpha_2} t_{\alpha_1} t_{\gamma}$, 
              where $\beta$ and $\gamma$ are simple closed curves in
              the diagrams \ref{Thm1-2ProofPALF_A} and  \ref{Thm1-2ProofPALF_B}.
    \end{enumerate}    
\end{theorem}

\section{Preliminaries}

    \subsection{Mapping class groups}
\label{sec:MCG}

In this subsection, we review a precise definition of the 
mapping class groups of surfaces with boundary and that of
Dehn twists along simple closed curves on surfaces.

\begin{definition}
    Let $F$ be a compact oriented connected surface with boundary.
    Let Diff${}^+(F, \partial F)$ be
    the group of all orientation-preserving self-diffeomorphisms
    of $F$ fixing the boundary $\partial F$ point-wise.
    Let Diff${}^+_0(F, \partial F)$ be
    the subgroup of Diff${}^+(F, \partial F)$ consisting of
    self-diffeomorphisms isotopic to the identity.
    The quotient group
    Diff${}^+(F, \partial F)/$ Diff${}^+_0(F, \partial F)$
    is called the mapping class group of $F$ and it is
    denoted by Map$(F, \partial F)$.
\end{definition}

\begin{definition}
    A {\em positive $($or right-handed$)$ Dehn twist}
    along a simple closed curve $\alpha$, $t_\alpha:F \rightarrow F$
    is a diffeomorphism  obtained by
    cutting $F$ along $\alpha$, twisting $360^\circ$
    to the right and regluing.
\end{definition}

  \subsection{PALF}

\begin{definition}\label{aLFdef}
    Let $M^4$ and $B^2$ be compact oriented smooth manifolds
    of dimensions $4$ and $2$.
    Let $f: M \rightarrow B$ be a smooth map.
    $f$ is called a {\em positive Lefschetz fibration} over $B$
    if it satisfies the following conditions (1) and (2):

    \begin{enumerate}
	\renewcommand{\labelenumi}{(\arabic{enumi})}
	\setlength{\itemsep}{0cm}
			\setlength{\parskip}{1mm}
	\item There are finitely many critical values
			$b_1, \ldots, b_m $ of $f$ in the interior of $B$
			and there is a unique critical point $p_i$
			on each fiber $f^{-1}(b_i)$, and
	\item The map $f$ is locally written as
			$f(z_1, z_2) = z_1^2 + z_2^2$ with respect to
			some local complex coordinates around
			$p_i$ and $b_i$ compatible with the orientations of
			$M$ and $B$.
    \end{enumerate}
\end{definition}

\begin{definition}
    A positive Lefschetz fibration is called {\em allowable} if
    its all vanishing cycles are homologically non-trivial on the
    fiber.
    A positive allowable Lefschetz fibration over $D^2$ with
    bounded fibers is called a {\em PALF} for short.
\end{definition}

The following Lemma is useful to prove Theorem \ref{MainThm:plug1}.

\begin{lemma}[{cf. Akbulut-Ozbagci \cite[Remark 1]{AO}}]
    \label{Lefschetz_2_h_PALF}
    Suppose that a $4$-manifold $X$ admits a PALF.
    If a $4$-manifold $Y$ is obtained from $X$ by attaching a
    Lefschetz $2$-handle, then $Y$ also admits a PALF.
\end{lemma}

The Lefschetz $2$-handle is defined as follows.

\begin{definition}
    Suppose that $X$ admits a PALF.
    A {\em Lefschetz $2$-handle} is a $2$-handle attached
    along a homologically non-trivial simple closed curve
    in the boundary of $X$ with framing $-1$ relative to
    the product framing induced by the fiber structure.
\end{definition}

   \subsection{Stein surfaces}

In this subsection, we recall a definition of the Stein surfaces.
The question of which smooth $4$-manifolds admit
Stein structures can be completely reduced to a
problem in handlebody theory.

\begin{definition}
    A complex manifold is called a {\em Stein manifold} if
    it admits a proper biholomorphic embedding to
    $ \mathbb{C}^n $.
\end{definition}

\begin{definition}
    Let $W$ be a compact manifold with boundary.
    The manifold $W$ is called a {\em Stein domain} if it
    satisfies following condition:
    There is a Stein manifold $X$ and
    a plurisubharmonic function
    $\varphi : X \rightarrow [0, \infty)$
    such that
    $W = \varphi^{-1}([0, a])$
    for a regular value $a$ of $\varphi$.
\end{definition}

\begin{definition}
    A Stein manifold or a Stein domain is called
    a {\em Stein surface} if its complex dimension is $2$.
\end{definition}

\subsection{Plugs}

In this subsection, we give the definition of the plug.

\begin{definition}(Akbulut-Yasui \cite[Definition 2.2.]{AY1})
    Let $P$ be a compact Stein 4-manifold with boundary and
    $\tau : \partial P \rightarrow \partial P$ 
    an involution on the boundary, 
    which cannot extend to any self-homeomorphism of $P$. 
    We call $(P, \tau)$ a {\em Plug} of $X$, 
    if $P \subset X$ and $X$ keeps its homeomorphism type and changes its
    diffeomorphism type when removing $P$ and gluing it via $\tau$. 
    We call $(P, \tau)$ a {\em Plug} if there exists a smooth
    4-manifold $X$ such that $(P, \tau)$ is a plug of $X$. 
\end{definition}

\begin{definition}(Akbulut-Yasui \cite[Definition 2.3.]{AY1})
    Let $W_{m, n}$ be a smooth 4-manifold
    given by Figure \ref{AY_plug_KirbyDiam}.
    Let $f_{m, n}: \partial W_{m, n} \rightarrow \partial W_{m, n}$
    be the obvious involution obtained from first surgering
    $S^1 \times D^3$ to $D^2 \times S^2$ in the
    interiors of $W_{m, n}$, then surgering the
    other imbedded $D^2 \times S^2$ back to
    $S^1 \times D^2$ (i.e. replacing the dot in Figure
    \ref{AY_plug_KirbyDiam}). 
\end{definition}

\begin{theorem}(Akbulut-Yasui \cite[Theorem 2.5(2)]{AY1})
    For $m \geq 1$ and $n \geq 2$, the pair $(W_{m, n}, f_{m, n})$ is
    a plug.
\end{theorem}

\section{Proofs of Theorems \ref{MainThm:plug1} and \ref{MainThm:plug2}. }

In this section, we give the 
proof of  Theorem \ref{MainThm:plug1} and Theorem \ref{MainThm:plug2}. \\

\begin{figure}
    \begin{center}
	\includegraphics[width=129mm]{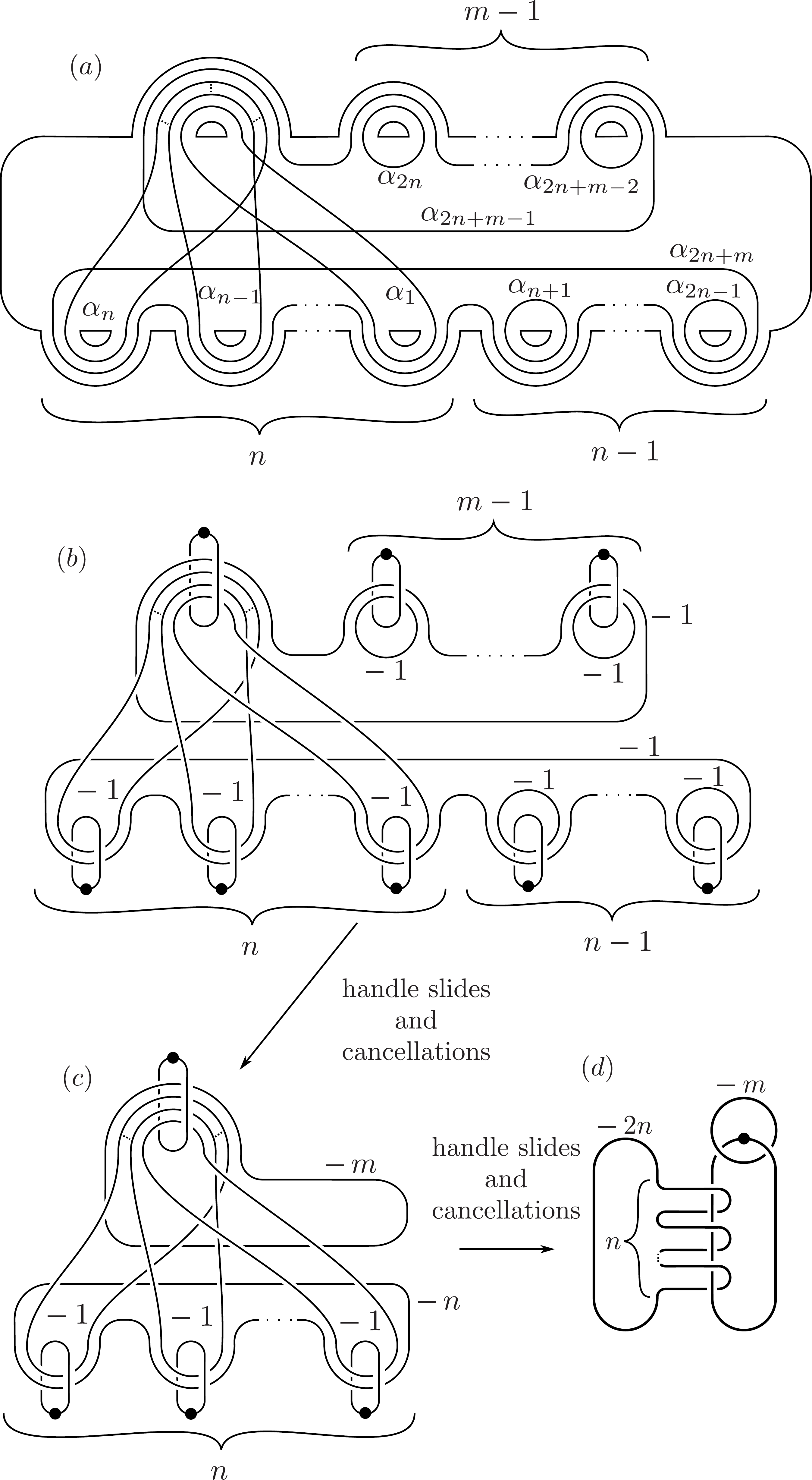}
    \end{center}
    \caption{}
    \label{Proof_plug2}
\end{figure}

{\it Proof of Theorem  \ref{MainThm:plug1}.}
Let $F_{m,n}$ be the compact oriented surface of genus zero
with $2n+m$ boundary components 
and $\alpha_1,\ldots ,\alpha_{2n+m}$ the curves on 
$F_{m,n}$ shown in Figure \ref{Proof_plug2} $(a)$. 
Note that Figure \ref{AkYsPlug} and 
Figure \ref{Proof_plug2} $(a)$ show the same PALF.
We denote the right-handed Dehn twists along $\alpha_1,\ldots ,\alpha_{2n+m}$ 
by $t_{\alpha_1}, \ldots ,t_{\alpha_{2n+m}}$, respectively. 
Let $f:X_{m,n}\rightarrow D^2$ be a Lefschetz fibration over $D^2$ 
with monodromy representation $(t_{\alpha_{2n+m}},\ldots ,t_{\alpha_1})$. 
Since each curve $\alpha_i$ is homologically non-trivial on $F_{m,n}$, 
we see that $f$ is a PALF with fiber $F_{m,n}$. 

We now show that $X_{m,n}$ is diffeomorphic to $W_{m,n}$.

The Kirby diagram for $X_{m,n}$ corresponding to the 
monodromy representation $(t_{\alpha_{2n+m}},\ldots ,t_{\alpha_1})$ 
is given by Figure \ref{Proof_plug2} $(b)$. 
We slide the $-1$-framed $2$-handles over $-1$-framed $2$-handles
and erase canceling $1$-handle/$2$-handle pairs
to get Figure \ref{Proof_plug2} $(c)$.
We get Figure \ref{Proof_plug2} $(d)$ 
by sliding the $-m$-framed $2$-handle over $-1$-framed $2$-handles and
sliding the $-n$-framed $2$-handle over $-1$-framed $2$-handles
and erasing canceling $1$-handle/$2$-handle pairs.

\begin{figure}
    \begin{center}
	\includegraphics[width=80mm]{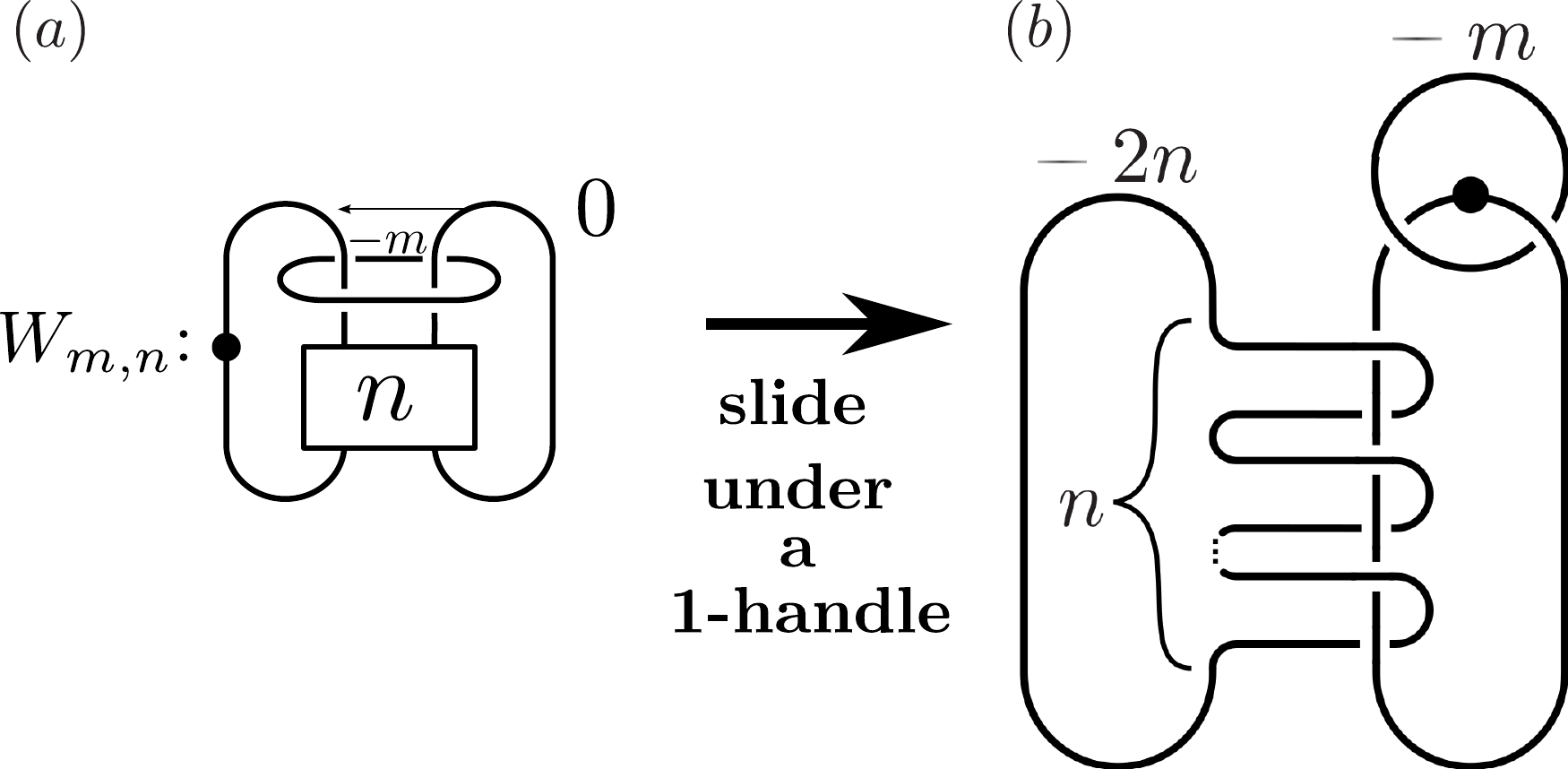}
    \end{center}
    \caption{}
    \label{Proof_plug1}
\end{figure}

The Kirby diagram for $W_{m,n}$ is given by Figure \ref{Proof_plug1} $(a)$. 
We slide the $0$-framed $2$-handle under the $1$-handle to get 
Figure \ref{Proof_plug1} $(b)$. 

Since Figure \ref{Proof_plug1} $(b)$ and Figure \ref{Proof_plug2} $(d)$ are the same, 
we conclude that $X_{m,n}$ is diffeomorphic to $W_{m,n}$, 
which implies the theorem. 
$\square$ \\

\noindent
{\itshape Proof of Theorem \ref{MainThm:plug2}. } 
(1) The plug twist of $A$ along
$(W_{1, 2}, f_{1, 2})$ is represented by
replacing the dot with 0 mutually
by the definition of the Akbulut-Yasui plug.
Therefore $B$ is obtained from $A$ by
plug twisting along $(W_{1, 2}, f_{1, 2})$.

\begin{figure}[htbp]
    \centering
    \begin{overpic}[width=4.5cm]{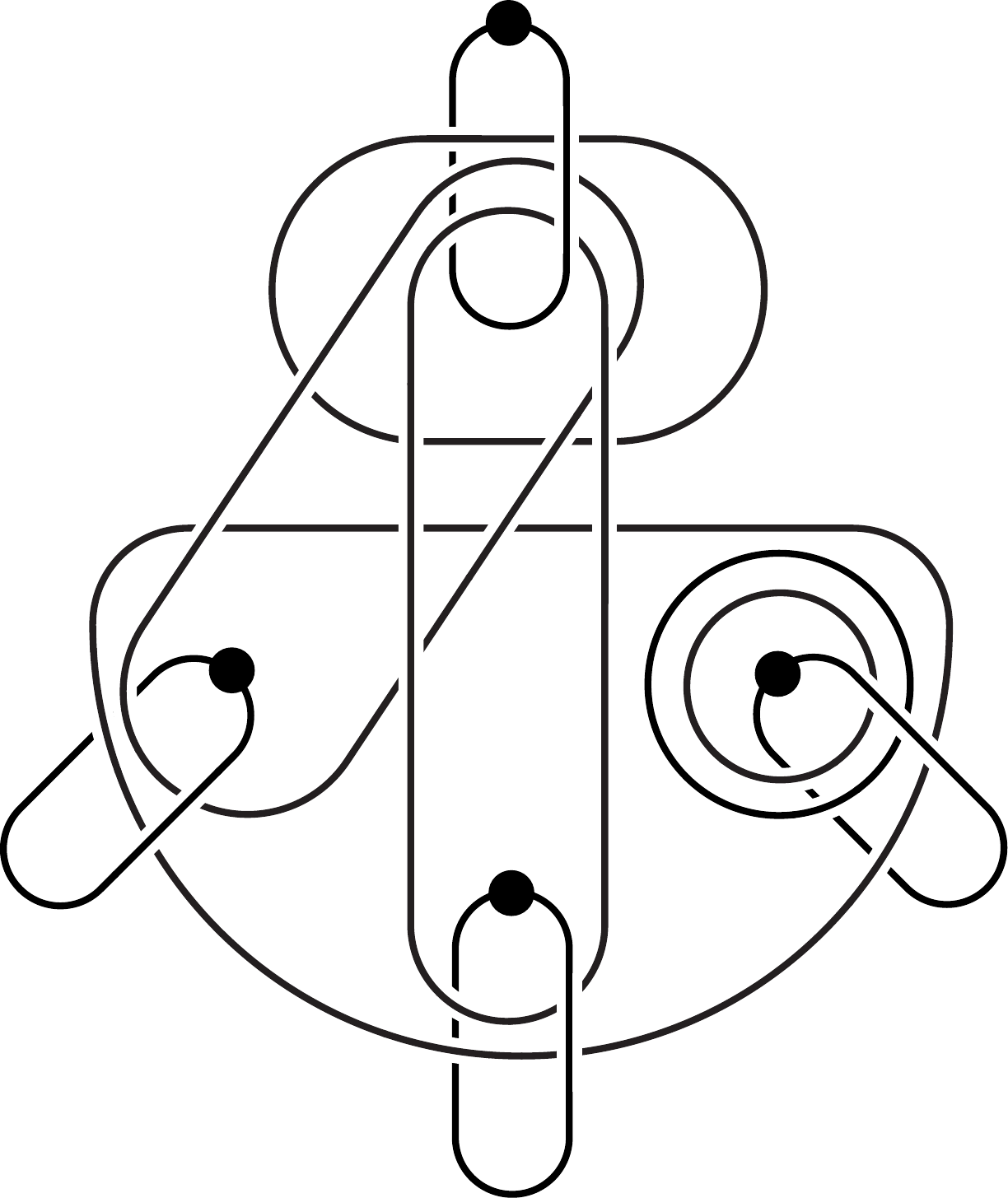}
	\put(16,85){$-1$}
	\put(10,62){$-1$}
	\put(12,12){$-1$}
	\put(39,35){$-1$}
	\put(56,25){{\small $-1$}}
	\put(60.5,46.5){{\tiny $-1$}}
    \end{overpic}
    \caption{}
    \label{diam:plugAKirby}
\end{figure}

\begin{figure}[htbp]
    \centering
    \begin{overpic}[width=12cm]{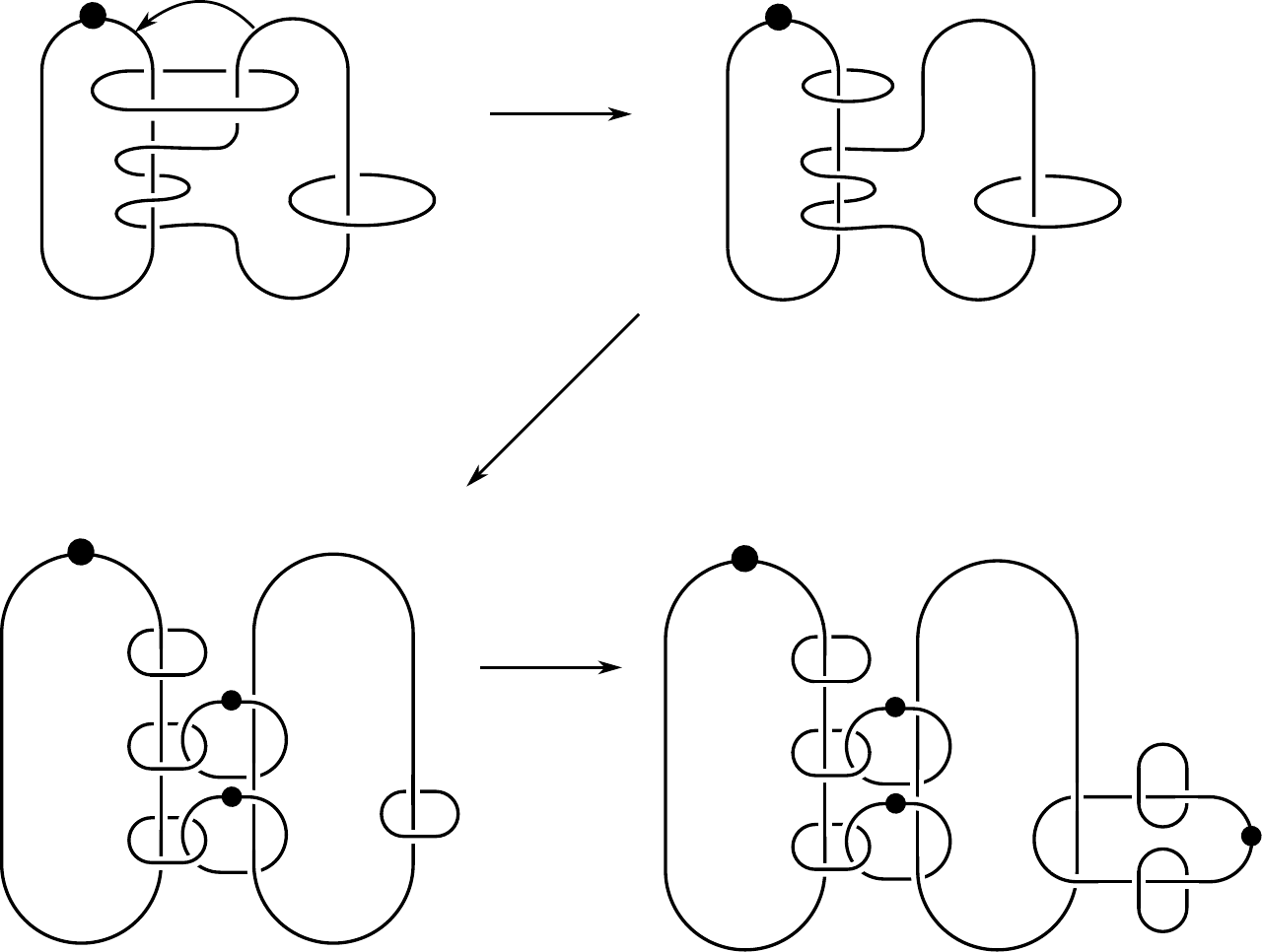}
	\put(12,76){\parbox{1cm}{slide} }
	\put(12.8,70.8){$-1$}
	\put(30,62){$-2$}
	\put(29,69){$0$}
	\put(40,68){slide}
	\put(67,71){$-1$}
	\put(84,69){$-4$}
	\put(84,62){$-2$}
	\put(30,46){\parbox{1.7cm}{creating canceling pairs}}
	\put(13,26.5){$-1$}
	\put(5,15.5){$-1$}
	\put(5,8){$-1$}
	\put(34,26.5){$-2$}
	\put(34,13){$-2$}
	\put(39,14){\parbox{1.3cm}{{\footnotesize creating a canceling pair}}}
	\put(86,26.5){$-1$}
	\put(65,26.5){$-1$}
	\put(57,15.5){$-1$}
	\put(57,8){$-1$}
	\put(90,18){$-1$}
	\put(90,-1){$-1$}
	\put(0,75){(a)}
	\put(51,75){(b)}
	\put(0,33){(c)}
	\put(51,33){(d)}
    \end{overpic}
    \caption{}
    \label{diam:AProof}
\end{figure}

\noindent
(2) First, we transform the Kirby diagram of $A$
as in Figure~\ref{diam:AProof}. 
Figure~\ref{diam:AProof} $(a)$ is the Kirby diagram of $A$. 
We slide the $0$-framed $2$-handle under the $1$-handle to get 
Figure~\ref{diam:AProof} $(b)$.
We get Figure~\ref{diam:AProof} $(c)$ by creating canceling 
$1$-handle/$2$-handle pairs.
We create canceling pairs to get Figure~\ref{diam:AProof} $(d)$.
Then we consider the $4$-manifold with a genus zero PALF structure
as in Figure~\ref{Thm1-2ProofPALF_A}. 
The obvious Kirby diagram for this manifold is given by 
Figure~\ref{diam:plugAKirby}. 
Therefore, the manifold $A$ admits a genus zero PALF structure.

\begin{figure}[htbp]
    \centering
    \begin{overpic}[width=4.5cm]{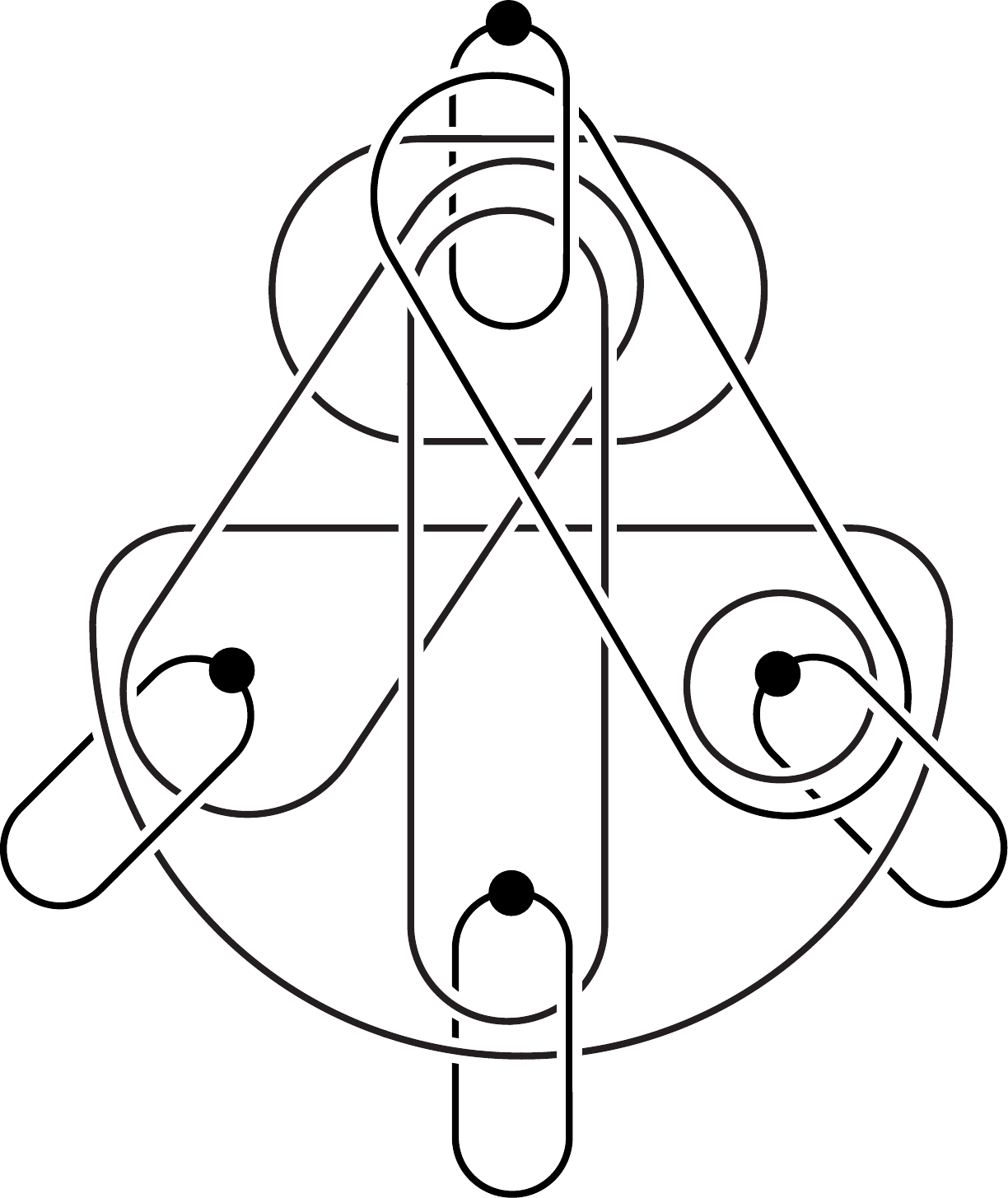}
	\put(12,85){$-1$}
	\put(10,62){$-1$}
	\put(12,12){$-1$}
	\put(39,35){$-1$}
	\put(54,50.5){{\small $-1$}}
	\put(68,62){$-1$}
    \end{overpic}
    \caption{}
    \label{diam:plugBKirby}
\end{figure}

\begin{figure}[htbp]
    \centering
    \begin{overpic}[width=12cm]{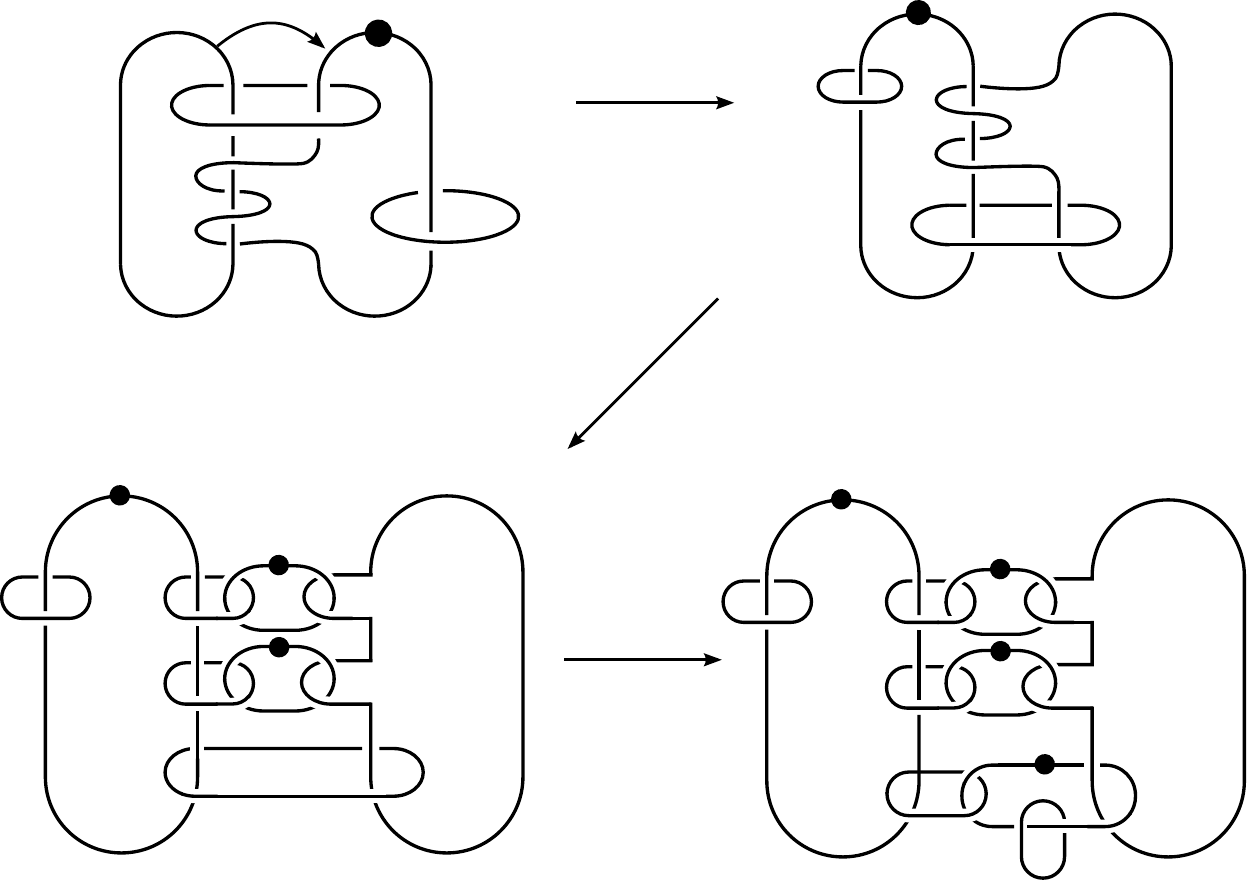}
	\put(10,69){$0$}
	\put(19,65){$-1$}
	\put(36,56){$-2$}
	\put(48,64){\parbox{1cm}{slide} }
	\put(62,65.5){$-1$}
	\put(95,65){$-4$}
	\put(79.5,47){$-2$}
	\put(35,40){\parbox{1.7cm}{creating canceling pairs} }
	\put(-2,25){$-1$}
	\put(10,25){$-1$}
	\put(10,18){$-1$}
	\put(10,10.5){$-2$}
	\put(43,25){$-2$}
	\put(44,10){\parbox{1.6cm}{creating a canceling pair} }
	\put(55,25){$-1$}
	\put(66.2,23){$-1$}
	\put(66.2,16){$-1$}
	\put(66.2,8){$-1$}
	\put(77,0){$-1$}
	\put(92,31){$-1$}
	\put(3,69){(a)}
	\put(58,69){(b)}
	\put(0,31){(c)}
	\put(58,31){(d)}
    \end{overpic}
    \caption{}
    \label{diam:BProof}
\end{figure}

Similarly, we transform the Kirby diagram of $B$
as in Figure~\ref{diam:BProof}.
Figure~\ref{diam:BProof} $(a)$ is the Kirby diagram of $B$. 
We slide the $0$-framed $2$-handle under the $1$-handle to get 
Figure~\ref{diam:BProof} $(b)$.
We get Figure~\ref{diam:BProof} $(c)$ by creating canceling 
$1$-handle/$2$-handle pairs.
We create canceling pairs to get Figure~\ref{diam:BProof} $(d)$.
Then we consider a $4$-manifold which admits a 
genus zero PALF structure as in Figure~\ref{Thm1-2ProofPALF_B}.
The obvious Kirby diagram for this manifold is given by 
Figure~\ref{diam:plugBKirby}. 
Therefore, the manifold $A$ admits a genus zero PALF structure.

\noindent
(3)We give a handle decomposition of $A$ and $B$ with
one $0$-handle, two $2$-handles as in
Figure~\ref{MfdAB}.
Therefore, 
$H_0(A;\mathbb{Z}) \cong H_0(B;\mathbb{Z}) \cong \mathbb{Z}$,
$H_1(A;\mathbb{Z}) \cong H_1(B;\mathbb{Z}) \cong \{0\}$,
$H_2(A;\mathbb{Z}) \cong H_2(B;\mathbb{Z}) \cong
\mathbb{Z}\oplus\mathbb{Z}$, and
$H_i(A;\mathbb{Z}) \cong H_i(B;\mathbb{Z}) \cong \{0\}\  (i \geq 3)$. 
The second Betti numbers of $A$ and $B$ are equal to 2. 

\noindent
(4) By the Kirby diagrams of $A$ and $B$ (Figure~\ref{MfdAB}),
both of the  boundaries
of $A$ and $B$ are represented by integral surgery
diagrams.
Therefore the boundaries of $A$ and $B$ are diffeomorphic to each other. 

\begin{figure}[htbp]
    \centering
    \begin{overpic}[width=10cm]{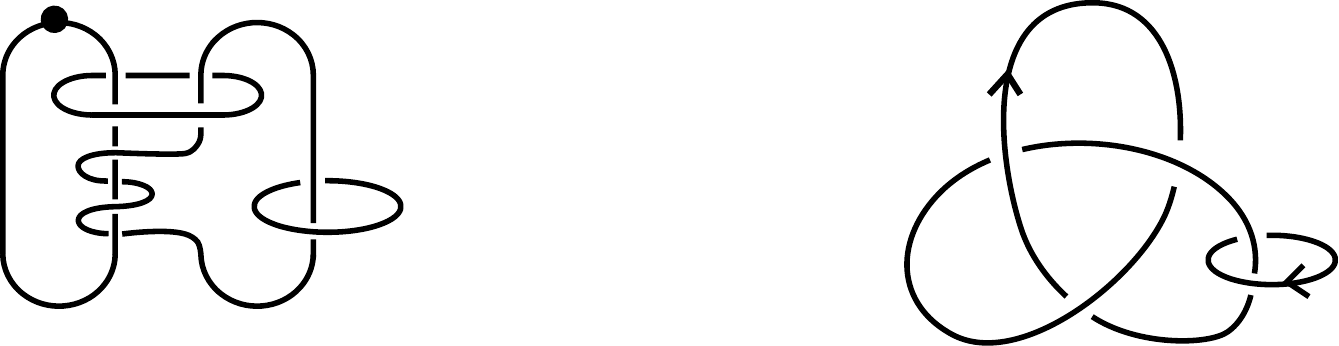}
	\put(10,21){$-1$}
	\put(25,21){$0$}
	\put(25,13){$-2$}
	\put(10,1){$A$}
	\put(37,18){\vector(1,0){20}}
	\put(38,8){\parbox{1.9cm}{handle slides and cancellation} }
	\put(68,20){$-8$}
	\put(94,9){$-2$}
    \end{overpic}
    \caption{}
    \label{diam:thm104:3:proof01}
\end{figure}

\begin{figure}[htbp]
    \centering
    \begin{overpic}[width=11cm]{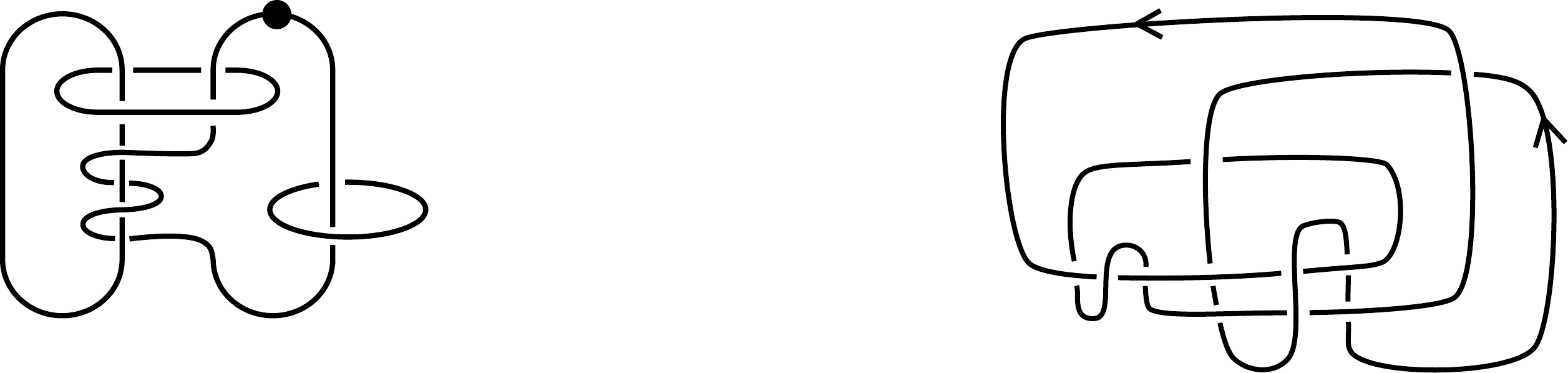}
	\put(9,21){$-1$}
	\put(0,24){$0$}
	\put(25,13){$-2$}
	\put(10,1){$B$}
	\put(37,18){\vector(1,0){20}}
	\put(38,8){\parbox{1.9cm}{handle slides and cancellation} }
	\put(66,18){$-8$}
	\put(82,16){$-3$}
    \end{overpic}
    \caption{}
    \label{diam:thm104:3:proof02}
\end{figure}

\noindent
(5)
We transform the Kirby diagrams of the manifolds $A$ and $B$
by Kirby calculus in 
Figure~\ref{diam:thm104:3:proof01} and Figure~\ref{diam:thm104:3:proof02}.
We obtain the intersection matrices
\begin{displaymath}
     \left(
    \begin{array}{cc}
        -8 & 1 \\
        1 & -2 \\
    \end{array}
    \right)
and
        \left(
    \begin{array}{cc}
        -8 & -3 \\
        -3 & -3 \\
    \end{array}
    \right)
\end{displaymath}
of $A$ and $B$ from the diagrams, respectively.
The former is even and the latter is odd.
Therefore $A$ and $B$ do not have isomorphic intersection form,
especially $A$ and $B$ are not homeomorphic. 

\noindent
(6)
The genus zero PALF structure on $W_{1, 2}$
is obtained from a trivial surface bundle over $D^2$
by attaching Lefschetz $2$-handles along simple closed
curves in Figure~\ref{W12}.
The genus zero PALF structure  $A$ (respectively $B$)
is obtained from surface bundle over $D^2$ by
attaching Lefschetz $2$-handles along simple closed
curves in Figure~\ref{Thm1-2ProofPALF_A}
(respectively Figure~\ref{Thm1-2ProofPALF_B}).
Therefore the monodromy representation of the
PALF on $A$ (respectively $B$) is 
$t_{\alpha_5} t_{\alpha_4} t_{\alpha_3} t_{\alpha_2} t_{\alpha_1} t_{\beta}$
(respectively $t_{\alpha_5} t_{\alpha_4} t_{\alpha_3} t_{\alpha_2}  t_{\alpha_1} t_{\gamma}$)
(where $t_{\alpha_i}$ is right-handed Dehn twist along
the simple closed curve $\alpha_i$).
\hfill $\Box$

\begin{figure}[htbp]
    \centering
    \begin{overpic}[width=4.5cm]{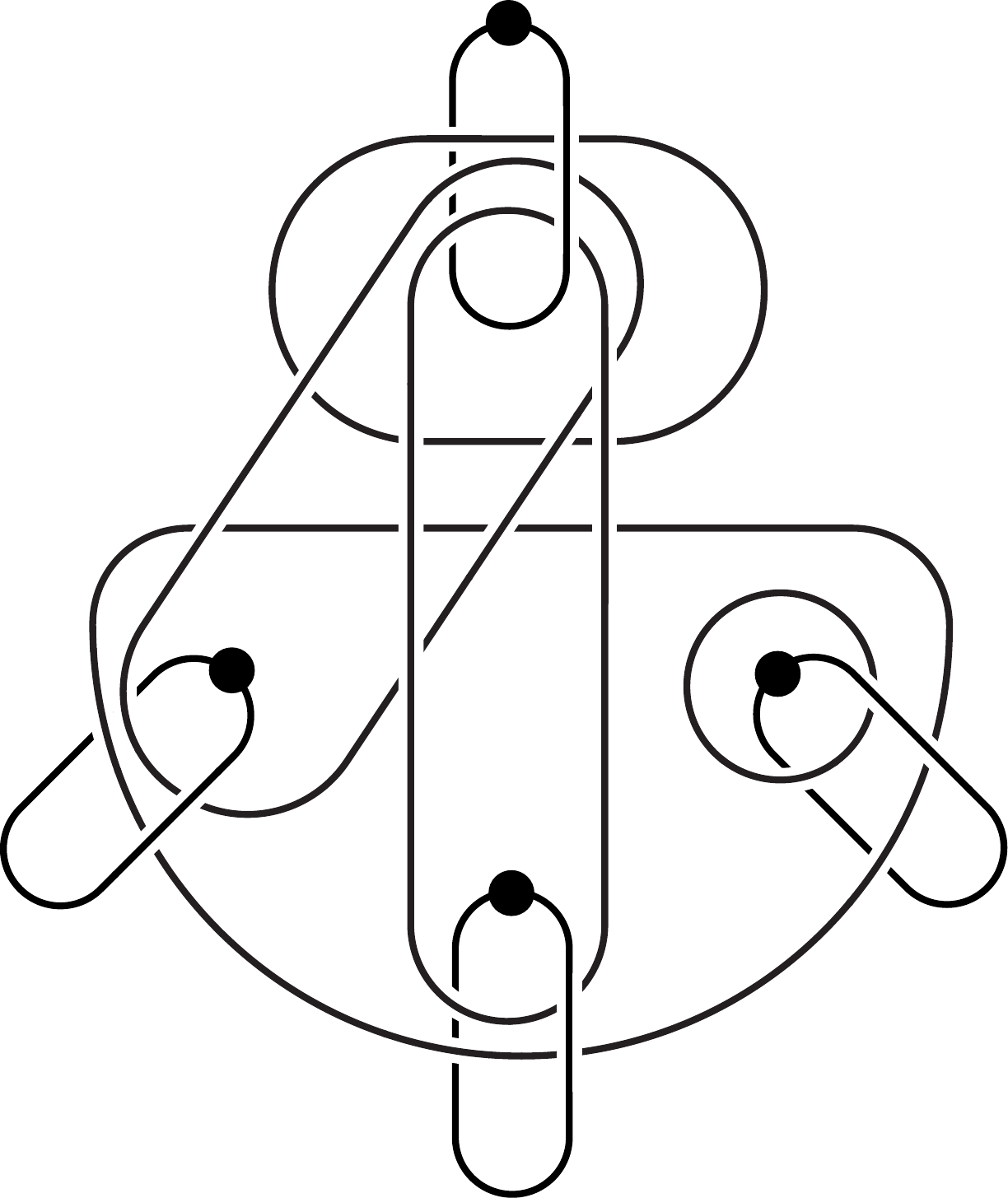}
	\put(12,85){$-1$}
	\put(10,62){$-1$}
	\put(12,12){$-1$}
	\put(39,35){$-1$}
	\put(56,29){{\small $-1$}}
    \end{overpic}
    \caption{}
    \label{diam:plugW12Kirby}
\end{figure}

\begin{figure}[htbp]
    \centering
    \begin{overpic}[width=3cm]{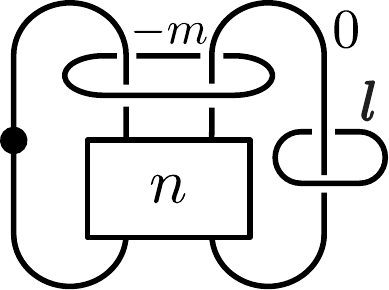}
    \end{overpic}
    \caption{}
    \label{diam16Wmnl}
\end{figure}

\begin{remark}
    Theorem \ref{MainThm:plug2} $(2)$ can be generalized in the same way for any 
    $m \geq 1, n \geq 2, l \leq -2$ in Figure \ref{diam16Wmnl}.
\end{remark}

\clearpage

\begin{ack}
    The author would like to thank his adviser
    Hisaaki Endo for his helpful comments and
    his encouragement.
    The author wishes to thank Kouichi Yasui 
    and Yuichi Yamada for
    their useful comments.
\end{ack}

\end{document}